\documentclass[a4paper,10pt,reqno]{amsart}
\usepackage[a4paper,tmargin=40mm,bmargin=35mm,hmargin=30mm]{geometry}
\usepackage[T1]{fontenc}
\usepackage{lmodern}
\usepackage{amsmath}
\usepackage{amssymb}
\usepackage{amsthm}
\usepackage{stmaryrd}
\usepackage{tikz}
\usepackage{booktabs}
\usepackage{multirow}
\usepackage{multicol}
\usepackage{here}
\usepackage{aliascnt}
\usepackage{hyperref}
\usepackage[noabbrev]{cleveref}
\usepackage[all]{xy}
\newcommand{\NewTheorem}[2]{
	\newaliascnt{#1}{TheoremEnvironment}
	\newtheorem{#1}[#1]{#1}
	\aliascntresetthe{#1}
	\crefname{#1}{#1}{#2}
	\Crefname{#1}{#1}{#2}
}

\theoremstyle{definition}

\NewTheorem{Definition}{Definitions}
\theoremstyle{remark}
\NewTheorem{Remark}{Remarks}
\NewTheorem{Axiom}{Axioms}
\NewTheorem{Example}{Examples}
\NewTheorem{Observation}{Observations}
\NewTheorem{Convention}{Conventions}
\NewTheorem{Remark and Notation}{Remark and Notations}
\NewTheorem{Setting}{Settings}
\NewTheorem{Question}{Questions}
\NewTheorem{Answer}{Answers}
\NewTheorem{Conjecture}{Conjectures}
\NewTheorem{Problem}{Problems}
\NewTheorem{Solution}{Solutions}
\NewTheorem{Goal}{Goals}
\NewTheorem{Comment}{Comments}
\NewTheorem{Aim}{Aims}
\NewTheorem{Caution}{Cautions}
\NewTheorem{Exercise}{Exercises}

\theoremstyle{plain}
\NewTheorem{Proposition}{Propositions}
\NewTheorem{Lemma}{Lemmas}
\NewTheorem{Theorem}{Theorems}
\NewTheorem{Corollary}{Corollaries}

\crefname{enumi}{}{}
\Crefname{enumi}{}{}
\creflabelformat{enumi}{(#2#1#3)}
\crefname{enumii}{}{}
\Crefname{enumii}{}{}
\creflabelformat{enumii}{(#2#1#3)}
\crefname{enumiii}{}{}
\Crefname{enumiii}{}{}
\creflabelformat{enumiii}{(#2#1#3)}

\makeatletter
\renewcommand{\p@enumii}{}
\renewcommand{\p@enumiii}{}
\makeatother

\numberwithin{equation}{section}
\crefname{equation}{}{}
\Crefname{equation}{}{}
\creflabelformat{equation}{(#2#1#3)}

\newcommand{\SwapSymbols}[1]{
	\expandafter\let\expandafter\temporarysymbol\csname #1\endcsname
	\expandafter\let\csname #1\expandafter\endcsname\csname var#1\endcsname
	\expandafter\let\csname var#1\endcsname\temporarysymbol
}
\usepackage{hyperref}
\hypersetup{%
  bookmarksnumbered=true,%
  colorlinks=true,%
  linkcolor=blue,%
  citecolor=blue,%
  filecolor=blue,%
  menucolor=blue,%
  urlcolor=blue,%
  bookmarksopen=true,%
  bookmarksdepth=2,%
  pageanchor=true}
\usepackage[noabbrev]{cleveref}

\SwapSymbols{epsilon}
\SwapSymbols{phi}
\SwapSymbols{Gamma}
\SwapSymbols{Delta}
\SwapSymbols{Theta}
\SwapSymbols{Lambda}
\SwapSymbols{Xi}
\SwapSymbols{Pi}
\SwapSymbols{Sigma}
\SwapSymbols{Upsilon}
\SwapSymbols{Phi}
\SwapSymbols{Psi}
\SwapSymbols{Omega}

\newcommand{\cA}{\mathcal{A}}

\newcommand{\cE}{\mathcal{E}}
\newcommand{\cF}{\mathcal{F}}

\newcommand{\cI}{\mathcal{I}}

\newcommand{\cO}{\mathcal{O}}

\newcommand{\cS}{\mathcal{S}}

\newcommand{\cU}{\mathcal{U}}
\newcommand{\cV}{\mathcal{V}}

\newcommand{\cX}{\mathcal{X}}

\newcommand{\kF}{\mathfrak{F}}

\newcommand{\To}{\rightarrow}

\DeclareMathOperator{\Hom}{Hom}

\DeclareMathOperator{\Zg}{Zg}
\DeclareMathOperator{\fg}{fg-}

\DeclareMathOperator{\Mod}{Mod}

\DeclareMathOperator{\ZSupp}{\cO}
\DeclareMathOperator{\ZASpec}{ZASpec}

\DeclareMathOperator{\GKdim}{KG-dim}
\DeclareMathOperator{\fp}{fp-}

\DeclareMathOperator{\Ann}{Ann}

\DeclareMathOperator{\Ker}{Ker}

\let\Im\relax
\DeclareMathOperator{\Im}{Im}

\DeclareMathOperator{\Spec}{Spec}

\DeclareMathOperator{\ASpec}{ASpec}

\DeclareMathOperator{\AAss}{AAss}
\DeclareMathOperator{\ASupp}{ASupp}

\DeclareMathOperator{\Sp}{Sp}

\title{Classifying localizing subcategories of a locally coherent category}


\author{Reza Sazeedeh}

\address{Department of Mathematics, Urmia University, P.O.Box: 165, Urmia, Iran}
\email{rsazeedeh@ipm.ir}

\begin{document}

\begin{abstract}
 Let $\cA$ be a locally coherent Grothendieck category, $\fp\cA$ be the full subcategory of $\cA$ consisting of finitely presented objects  and $\ASpec\cA$ be the atom spectrum of $\cA$. In this paper, we classify localizing subcategories of finite type of $\cA$ via open subsets a new topology on $\ASpec\cA$. We investigate $\ASpec\fp\cA$ and show that if $\ASpec\fp\cA=\ASpec\cA$, then $\cA$ is locally noetherian. As an application, we specialize our investigation to the case of commutative coherent rings.
\end{abstract}

\maketitle
\tableofcontents

\section{Introduction}
 Throughout this paper, let $\cA$ be a Grothendieck category.  The idea of classifying subcategories of $\cA$ using a spectrum originates from the classification of subcategories of modules over a commutative ring $A$ via $\Spec A$, as established by numerous authors (\cite{Hop, Hov, GP, T}). Gabriel \cite{G} initiated the extension of this approach by defining the spectrum $\Sp\cA$ as the set of  isomorphism classes of indecomposable injective objects. For commutative noetherian rings, these objects corresponds bijectively to prime ideals. 
 
An alternative topology on $\Sp\cA$ is that of the Ziegler spectrum $\Zg\cA$, originally defined by Ziegler \cite{Z} associating to a ring $A$, a topological space whose points are the isomorphism classes of indecomposable pure-injective $A$-modules. In a locally coherent category $\cA$ with subcategory $\fp\cA$ of finitely presented objects, Herzog \cite{H} classified the localizing subcategories of finite type of $\cA$ and the Serre subcategories of $\fp\cA$ via open subsets of $\Zg\cA$. Subsequently, Krause [K] provided  alternative classifications of these same subcategories.

Alongside these developments,  Kanda \cite{K} recently introduced the atom spectrum $\ASpec\cA$, for an abelian category $\cA$  that need not have enough injectives, and endowed it with a topology. This construction was inspired by Storrer's \cite{Sto} work on  monoform modules and their equivalence relation over non-commutative rings. While the spectrum has been extensively studied over locally noetherian categories (see \cite{K}), in this paper we investigate $\ASpec\cA$ within locally coherent categories. In particular, we provide a new classification of both the localizing subcategories of finite type of $\cA$ and the  Serre subcategories of $\fp\cA$, in terms of open subsets of a newly defined topology on $\ASpec \cA$. We then specialize our results to the case of commutative coherent rings.

 In Section 2, we study locally monoform categories. An abelian category $\cA$ is said to be locally monoform if every non-zero object in $\cA$ contains a monoform subobject. Unfortunately, Grothendieck categories are not locally monoform in general, which limits our ability to find out further insights about $\cA$ via $\ASpec\cA$. A well-known example is  locally noetherian categories; see \cite{K}. We show that if $\cA$ is a locally monoform category, then  there is a bijective correspondence between  localizing subcategories of $\cA$ open subsets of $\ASpec\cA$; see \cref{lmsem}. Furthermore, we also show that semi-noetherian categories are locally monoform; see \cref{crit}.

In Section 3, we study locally coherent categories. A Grothendieck category $\cA$ is said to be locally coherent if $\cA$ has a generating set of finitely presented objects and the full subcategory $\fp\cA$ of finitely presented objects in $\cA$ is abelian. For an object $M$ in $\cA$, the atom support $\ASupp M$ can be defined analogously to the support of a module over a commutative ring.  In Section 3, we assume that $\cA$ is locally coherent and we define a new topology on $\ASpec\cA$ in which $\{\ASupp M|\hspace{0.1cm} M\in\fp\cA\}$ forms a basis of open subsets of $\ASpec\cA$. We use the symbol $\ZASpec\cA$ instead of $\ASpec\cA$ with this topology. For any subcategory $\cX$ of $\fp\cA$ and any subset $\cU$ of $\ASpec \cA$, assume that $\ASupp\cX=\bigcup_{M\in\cX}\ASupp M$ and $\ASupp^{-1}\cU=\{M\in\cA|\hspace{0.1cm} \ASupp M\subseteq \cU\}$. The following theorem classifies Serre subcategories of $\fp\cA$ using open subsets of $\ZASpec\cA$.
   
\begin{Theorem}[\cref{corpropp}]
Let $\cA$ be a locally coherent category. Then the map $\cX\mapsto \ASupp\cX$  provides an inclusion-preserving bijective correspondence between Serre subcategories of $\fp\cA$ and open subsets of $\ZASpec\cA$. The inverse map is given by $\cU\mapsto \ASupp^{-1}\cU\cap\fp\cA$.
\end{Theorem}
As a conclusion, we prove that there is a bijective correspondence between open subsets of $\ZASpec\cA$ and localizing subcategories of finite type of $\cA$; see \cref{ccooo}. As $\fp\cA$ is abelian, $\ASpec \fp\cA$ can be  investigated independently. It is known as Cohen's theorem that if the prime ideals of a commutative ring $A$ are finitely generated, then $A$ is noetherian. In the following theorem we extend this theorem for locally coherent categories.

\begin{Theorem}[\cref{noeth}]
Let $\cA$ be a locally coherent category. If $\ASpec\fp\cA=\ASpec\cA$ (i.e. if every $\alpha\in\ASpec\cA$ has a monoform representative in $\fp\cA$), then $\cA$ is locally noetherian.
\end{Theorem}   
  In the following theorem, we provide a sufficient condition under which a localizing subcategory of $\cA$ is of finite type.
  
 \begin{Theorem}[\cref{fgftg}]
Let $\cA$ be a locally coherent category and let $\cX$ be a localizing subcategory of $\cA$ such that every atom in $\ASupp\cX$ has a  monoform representative in $\fp\cA$. Then $\cX$  is of finite type.
\end{Theorem}

For any object $M$ in $\cA$, set $\ZSupp(M)=\{I\in\Zg\cA|\hspace{0.1cm} \Hom(M,I)\neq 0\}.$  Herzog showed that over a locally coherent category, the set $\{\ZSupp(M)|\hspace{0.1cm} M\in\fp\cA\}$ can serve as a basis for a topology on $\cA$. This topology, denoted by $\Zg\cA$, is known as Ziegler topology of $\cA$. A category $\cA$ is said to be locally uniform if every non-zero object in $\cA$ has a uniform subobject. 

In Section 4, we show that $\cA$ is locally monoform if and only if $\cA$ is locally uniform and the map $\theta:\ZASpec\cA\To\Zg\cA$, given by  $\alpha\mapsto E(\alpha)$ is a homeomorphism. We provide an example which show that any locally coherent is not necessarily locally monoform; see \cref{eex}.
 Furthermore, we provide a sufficient condition under which a locally coherent category is semi-noetherian (and so a locally monoform); see \cref{semii}.  
  
In Section 5, we assume that $A$ is a commutative coherent ring. We denote $\ZASpec \Mod A$ and  $\Zg \Mod A$ by $\ZASpec A$ and $\Zg A$. For a spectral topological space $X$, Hochster \cite{Ho} endowed the underlying set with a new, dual topological by defining its open subsets as those of the form $Y=\bigcup_{i\in\Omega} Y_i$, where $X\setminus Y_i$ is a quasi-compact open subset of $X$ for each $i\in\Omega$. The symbol $X^\star$ denotes $X$  with the new topology. We write $\Spec^\star A$ for $(\Spec A)^\star$. We establish a bijective correspondence between open subsets of $\ZASpec^\star A$ and that of $\ZASpec A$; see \cref{pla}.  We use Gabriel topologies on $A$ to classify the localizing subcategories of finite type of $\Mod A$ over a commutative  coherent ring $A$. 

\begin{Theorem}[\cref{flfg}]
\begin{enumerate}
\item  For any finitely generated ideal $\frak a$, the assignments $$\ASupp A/\frak a\longleftrightarrow\Zg A\setminus \cO(A/\frak a)$$ induce an inclusion-reversing bijection between open subsets of $\ZASpec A$ and closed subsets of $\Zg A$.

\item  The assignment $D\stackrel{\kF}\mapsto \kF_D$ induces  an inclusion-reversing bijective  between closed subsets of $\Zg A$ and Gabriel topologies having bases of finitely generated ideals of $A$. The inverse map is given by $\kF\stackrel{D}\mapsto D(\kF)$.

\item  The assignment $D\stackrel{\kF}\mapsto \kF_D$ induces  an inclusion-reversing bijective  between closed subsets of $\Zg A$ and Gabriel topologies having bases of finitely generated ideals of $A$. The inverse map is given by $\kF\stackrel{D}\mapsto D(\kF)$.

\item  The assignment   $\kF\mapsto \cX_\kF=\overrightarrow{\sqrt{\langle A/\frak a|\hspace{0.1cm} \frak a\in \kF\cap\fg A \rangle}}$ induces  an inclusion-preserving bijective  between  Gabriel topologies having bases of finitely generated ideals of $A$ and localizing subcategories of finite type of $\Mod A$. The inverse map is given by $\cX\mapsto\kF_{\cX}=\{\frak a|\hspace{0.1cm}A/\frak a\in\cX\}$.
\end{enumerate}
 \end{Theorem}


\section{Categories with enough atoms}

 We begin this section with key definitions which are used throughout this paper.

\begin{Definition}
(1) An abelian category $\cA$ with a generator is said to be a {\it Grothendieck category} if it has arbitrary direct sums and direct limits of short exact sequence are exact, this means that if a direct system of short exact sequences in $\cA$ is given, then the induced sequence of direct limits is a short exact sequence. 

(2) An object $M$ of a Grothendieck category $\cA$ is {\it finitely generated} if whenever there are subobjects $M_i\leq M$ for $i\in I$ satisfying $M=\underset{i\in I}\Sigma M_i$, then there is a finite subset $J\subseteq I$ such that $M=\underset{i\in J}\Sigma M_i$.  A Grothendieck  category $\cA$ is said to be {\it locally finitely generated} if it has a small generating set of finitely generated objects.

(3) A Grothendieck category $\cA$ is said to be {\it locally noetherian} if it has a small generating set of noetherian objects.
\end{Definition}
 
 \medskip

\begin{Definition}
A full subcategory $\cX$ of an abelian category $\cA$ is said to be {\it Serre} if for any exact sequence $0\To M\To N\To K\To 0$ in $\cA$, the object $N$ belongs to $\cX$ if and only if $M$ and $K$ belong to $\cX$.
\end{Definition}

  \begin{Definition} 
A Serre subcategory $\cX$ of a Grothendieck category $\cA$ is said to be {\it localizing} if the
canonical functor $F:\cA\To\cA/\cX$ admits a right adjoint functor $G:\cA/\cX\To \cA$ (known as the {\it section functor}). The right adjoint of the inclusion functor $i:\cX\To\cA$, denoted by $t_\cX:\cA\To\cX$ is called {\it  the radical functor} associated to $\cX$. For any object $M$ in $\cA$, $t_\cX(M)$ is the largest subobject of $M$ contained in $\cX$.
\end{Definition}

\begin{Definition}
(i) A non-zero  object $U$ in $\cA$ is said to be {\it uniform} if every non-zero subobject of $U$ is an essential subobject of $U$.

 (ii) A non-zero object $M$ in $\cA$ is {\it monoform} if for any
non-zero subobject $N$ of $M$, there exists no common non-zero
subobject of $M$ and $M/N$ which means that there does not exist a
non-zero subobject of $M$ which is isomorphic to a subobject of
$M/N$. We denote by $\ASpec _0\cA$, the set of all monoform
objects in $\cA$.
\end{Definition}

 Two monoform objects $H$ and $H'$ in $\cA$ are said to be {\it
atom-equivalent} if they have a common non-zero subobject. The atom equivalence establishes an
equivalence relation on monoform objects. For every
monoform object $H$ in $\cA$,  we denoted by $\overline{H}$, the
 {\it equivalence class} of $H$, that is
\begin{center}
  $\overline{H}=\{G\in\ASpec _0\cA|\hspace{0.1cm} H \hspace{0.1cm}
 {\rm and}\hspace{0.1cm} G$ have a common non-zero subobject$\}.$
\end{center}

$\bullet$ In the rest of this paper, we always assume that $\cA$ is a Grothendieck category. 
   \begin{Definition}
 The {\it atom spectrum} $\ASpec\cA$ of $\cA$ is the quotient class
of $\ASpec _0\cA$ consisting of  equivalence classes induced by
this equivalence relation; in other words 
$$\ASpec\cA=\{\overline{H}|\hspace{0.1cm} H\in\ASpec _0\cA\}.$$
Any equivalence class is called an {\it
atom} of $\ASpec\cA$. We observe that since $\cA$ is a Grothendieck category, $\ASpec\cA$ is a set. 
\end{Definition}

\medskip
 The notion support and associated prime of  modules over a commutative ring can be generalized to objects in a Grothendieck category $\cA$ as follows:  

\begin{Definition}
\begin{enumerate}
 Let $M$ be an object in $\cA$.

\item The {\it atom support} of $M$, denoted by $\ASupp M$, is defined as 
$$\ASupp M=\{\alpha\in\ASpec\cA|\hspace{0.1cm} {\rm there\hspace{0.1cm} exists\hspace{0.1cm}}
 H\in\alpha\hspace{0.1cm} {\rm which \hspace{0.1cm}is\hspace{0.1cm} a
 \hspace{0.1cm}subquotient\hspace{0.1cm} of\hspace{0.1cm}} M\}.$$

\item The {\it associated atom} of $M$, denoted by $\AAss M$, is defined as
$$\AAss M=\{\alpha\in\ASupp M|\hspace{0.1cm} {\rm there\hspace{0.1cm} exists\hspace{0.1cm}}
 H\in\alpha\hspace{0.1cm} {\rm which \hspace{0.1cm}is\hspace{0.1cm} a
 \hspace{0.1cm}subobject\hspace{0.1cm} of\hspace{0.1cm}} M\}.$$ 
 
\item An atom $\alpha$ in $\ASpec\cA$ is said to be {\it maximal} if there
exists a simple object $H$ in $\cA$ such that $\alpha=\overline{H}$. The set of maximal atoms in $\ASpec\cA$ is denoted by $\textnormal{m-}\ASpec \cA$.
 
\item  A subset $\Phi$ of $\ASpec\cA$ is said to be {\it open} if for any
$\alpha\in\Phi$, there exists a monoform object $H$  in $\cA$ such that $\alpha=\overline{H}$ and $\ASupp H\subset\Phi$.
\end{enumerate}
\end{Definition}
\begin{Remark}\label{remop} For any object $M$ in $\cA$, it is clear that $\ASupp M$ is an open subset of $\ASpec\cA$. Then $\ASpec \cA$ is equipped with a topology having $\{\ASupp M|\hspace{0.1cm} M\in\cA\}$ as a basis of open subsets. For any subcategory $\cX$ of $\cA$, we set 
$$\ASupp\cX=\underset{M\in\cX}\bigcup\ASupp M$$ which is an open subset of $\ASpec\cA$. Also, for every subset $\cU$ of $\ASpec\cA$, we set $$\ASupp^{-1}\cU=\{M\in\cA|\hspace{0.1cm}\ASupp M\subseteq \cU\}$$ which is a Serre subcategory of $\cA$.  Furthermore, since $\cA$ is a Grothendieck category, it follows from \cite[Proposition 3.11]{K} that $\ASupp^{-1}\cU$ is a localizing subcategory of $\cA$.
 \end{Remark}

\medskip

 \begin{Definition}
 An object $M$ in $\cA$ is said to be {\it locally monoform} if every non-zero subquotient of $M$ contains a monoform subobject. An abelian category $\cA$ is said to be {\it locally monoform} if every non-zero object in $\cA$ is locally monoform; in other words any non-zero object in $\cA$ contains a monoform subobject.  
\end{Definition}

\medskip

 \begin{Lemma}\label{lemcru}
 Let $\cX$ be a localizing subcategory of $\cA$ and let $M$ be a locally monoform object in $\cA$. If $\ASupp M\subset\ASupp \cX$, then $M\in\cX$.
 \end{Lemma}
\begin{proof}
Assume that $M$ is not in $\cX$ and $t_\cX(M)$ is the largest subobject of $M$ belonging to $\cX$. Since $M$ is a locally monoform, $M/t_\cX(M)$  contains a monoform subobject $N/t_\cX(M)$. Therefore $\overline{N/t_\cX(M)}\in\ASupp\cX$ and so there exists an object $X\in\cX$ such that $\overline{N/t_\cX(M)}\in\ASupp X$. Thus $N/t_\cX(M)$ contains a non-zero subobject isomorphic to a subquotient of $X$. But this implies that $t_\cX(N/t_\cX(M))$ is non-zero which is a contradiction.  
\end{proof}

\medskip

\begin{Proposition}\label{lmsem}
Let $\cA$ be a locally monoform category. Then the map $\cX\mapsto \ASupp\cX$ provides an inclusion-preserving  bijective correspondence between  localizing subcategories of $\cA$ and  open subsets of $\ASpec\cA$. The inverse map is given by $\cU\mapsto \ASupp^{-1}\cU$.
\end{Proposition}
\begin{proof}
For a localizing subcategory $\cX$ of $\cA$ and an open subset $\cU$ of $\ASpec\cA$, it suffices to prove that $\ASupp^{-1}(\ASupp\cX)=\cX$ and $\ASupp(\ASupp^{-1}\cU)=\cU$.  But this follows easily from the definition and using \cref{lemcru}.
\end{proof}

It is clear that locally noetherian categories are locally monoform (see [K, Proposition 3.5]). In the rest of this section we show that semi-noetherian categories are another example of locally monoform categories.

\medskip
\begin{Definition}\label{Gab}
Let $\cA$ be a Grothendieck category.  We define the {\it Krull-Gabriel filtration} of $\cA$ as follows:
 
  For any ordinal (i.e ordinal number) $\sigma$
we denote by $\cA_{\sigma}$, the localizing subcategory of $\cA$ which is defined in the following manner:

$\cA_{-1}$ is the zero subcategory.

$\cA_0$ is the smallest localizing subcategory containing all simple objects.

Let us assume that $\sigma=\rho+1$ and denote by $F_{\rho}:\cA\To\cA/\cA_{\rho}$ the canonical functor and by 
$G_{\rho}:\cA/\cA_{\rho}\To\cA$ the right adjoint functor of $F_{\rho}$. Then an object $X$ in $\cA$ will belong to 
$\cA_{\sigma}$ if and only if $F_{\rho}(X)\in{\rm Ob}(\cA/\cA_{\rho})_0$. The left exact radical functor (torsion functor) corresponding to $\cA_{\rho}$ is denoted by $t_{\rho}$.
If $\sigma$ is a limit ordinal, then $\cA_{\sigma}$ is the localizing subcategory generated by all localizing subcategories $\cA_{\rho}$ with $\rho< \sigma$. It is clear that if $\sigma\leq \sigma'$, then $\cA_{\sigma}\subseteq\cA_{\sigma'}$. Moreover, since the class of all localizing subcategories of $\cA$ is a set, there exists an ordinal $\tau$ such that $\cA_{\sigma}=\cA_{\tau}$ for all $\sigma \geq \tau$. Let us put $\cA_{\tau}=\cup_{\sigma}\cA_{\sigma}$. Then $\cA$ is said to be {\it semi-noetherian}  if $\cA=\cA_{\tau}$. We also say that the localizing subcategories $\{\cA_{\sigma}\}_{\sigma}$ define the Krull-Gabriel filtration of $\cA$. We say that an object $M$ in $\cA$ has the {\it Krull-Gabriel dimension} defined or $M$ is {\it semi-noetherian} if $M\in{\rm Ob}(\cA_{\tau})$. The smallest ordinal  $\sigma$ so that $M\in{\rm Ob}(\cA_{\sigma})$ is denoted by $\GKdim M$. We observe that $\GKdim 0=-1$ and $\GKdim M\leq 0$ if and only if $\ASupp M\subseteq \textnormal{m-}\ASpec \cA$. 
\end{Definition}
We note that any locally noetherian Grothendieck category is semi-noetherian; see \cite[Chap. 5, Theorem 8.5]{Po}. To be more precise, If $\cA\neq\cA_{\tau}$, then $\cA/\cA_{\tau}$ is also locally noetherian and so it has a non-zero noetherian object $X$. Then $X$ has a maximal subobject $Y$ so that $S=X/Y$ is simple. Therefore, $\sigma=\overline{S}\in\ASupp(\cA/\cA_{\tau})_0$ which is a contradiction by the choice of $\tau$.

\medskip 
 
	\begin{Definition}
	Given an ordinal $\sigma\geq 0$, we recall from \cite{MR, GW} that a non-zero object $M$ in $\cA$ is  $\sigma$-{\it critical} provided $\GKdim M=\sigma$ while $\GKdim M/N<\sigma$ for all non-zero subobjects $N$ of $M$. It is clear that any non-zero subobject of a $\sigma$-critical object is $\sigma$-critical. An object $M$ is said to be {\it critical} if it is $\sigma$-critical for some ordinal $\sigma$. It is clear to see that any critical object is monoform.  
	\end{Definition}

\begin{Lemma}\label{nonlim}
Let $M$ be a $\sigma$-critical object in $\cA$. Then $\sigma$ is a non-limit ordinal. 
\end{Lemma}
\begin{proof}
Assume that $\sigma$ is a limit ordinal. By the definition $\cA_\sigma$ is generated by $\bigcup_{\rho<\sigma}\cA_\rho$ and since $\bigcup_{\rho<\sigma}\cA_\rho$ is closed under quotients it follows from [St, Chap VI, Proposition 2.5] that $M$ contains a non-zero subobject $N\in \bigcup_{\rho<\sigma}\cA_\rho$.  Then there exists some $\rho<\sigma$ such that $N\in\cA_\rho$ so that $N\subseteq t_{\rho}(M)$. This implies that $t_{\rho}(M)\neq 0$  and so $\GKdim t_{\rho}(M)\leq \rho$.  But $t_{\rho}(M)$ is $\sigma$-critical which is a contradiction.  
\end{proof}

\medskip

\begin{Lemma}\label{simcr}
Let $\sigma$ be a non-limit ordinal and let $M$ be an object in $\cA$. If $F_{\sigma-1}(M)$ is simple, then $M/t_{\sigma-1}(M)$ is $\sigma$-critical. 
\end{Lemma}
\begin{proof}
Observe that $F_{\sigma-1}(M)\cong F_{\sigma-1}(M/t_{\sigma-1}(M))$ and so we may assume that $t_{\sigma-1}(M)=0$. Let $N$ be a non-zero subobject of $M$. Then $F_{\sigma-1}(N)$ is non-zero because $t_{\sigma-1}(M)=0$. Since $F_{\sigma-1}(M)$ is simple, we have $F_{\sigma-1}(M/N)=0$ and hence $\GKdim M/N<\sigma$. On the other hand, by the definition and the fact that $F_{\sigma-1}(M)$ is simple, we have $\GKdim M=\sigma$. 
\end{proof}

\medskip

The following lemma is known for modules (for example see \cite[Lemma 15.8]{GW}); but since the dimension is slightly different, we reprove it.
\begin{Lemma}\label{excrit}
If $M$ is a non-zero object in $\cA$ with Krull-Gabriel dimension, then $M$ has a non-zero critical subobject (and so a non-zero monoform subobject)
\end{Lemma}
\begin{proof}
Since ordinals satisfy the descending chain condition, we can choose a non-zero subobject $N$ of $M$ of minimal Krull-Gabriel dimension $\sigma$. If $\sigma$ is limit ordinal, by 
[St, Chap VI, Proposition 2.5], $N$ contains a non-zero subobject in $\bigcup_{\rho<\sigma}\cA_\rho$ which contradicts the choice of $N$ with the minimal Krull-Gabriel dimension. Then $\sigma$ is non-limit ordinal and $t_{\sigma-1}(N)=0$. Since $F_{\sigma-1}(N)\in (\cA/\cA_{\sigma-1})_0$, it follows again from \cite[Chap VI, Proposition 2.5]{St} that $F_{\sigma-1}(N)$ contains a simple subobject $S$. Then 
$N$ contains a subobject $H$ such that $F_{\sigma-1}(H)=S$ by \cite[Chap 4, Corollary 3.10]{Po}. Now, \cref{simcr} implies that $H$ is a $\sigma$-critical.   
\end{proof}

\medskip

\begin{Corollary}\label{crit}
If $M$ is an object in $\cA$ with Krull-Gabriel dimension, then $M$ is loally monoform. In particular, if $\cA$ is a semi-noetherian category, then $\cA$ is locally monoform.
\end{Corollary}
\begin{proof}
If $M$ has Krull-Gabriel dimension, then every non-zero subquotient of $M$ has Krull-Gabriel dimension as well. Thus \cref{excrit} implies that every subquotient of $M$ contains a monoform subobject so that $M$ is locally monoform. The second assertion is clear. 
\end{proof}

\medskip
For every ordinal $\sigma$, the localizing subcategory $\cA_{\sigma}$ of $\cA$ is generated by critical objects.
\begin{Proposition}\label{critp}
Let $\sigma$ be an ordinal. Then $\cA_{\sigma}$ is generated by  all $\delta$-critical objects in $\cA$ with $\delta\leq \sigma$. 
\end{Proposition}
\begin{proof}
Let $\cX$ be the localizing subcategory of $\cA$ generated by all $\delta$-critical objects with $\delta\leq \sigma$. It is clear that $\cX\subseteq\cA_\sigma$. For the other direction, let $M\in\cA_\sigma$ and $t_\cX(M)$  be the largest subobject of $M$ belonging to $\cX$. If $M/t_\cX(M)$ is nonzero, by \cref{excrit}, it contains a critical subobject $N/t_\cX(M)$. Since $\GKdim M/t_\cX(M)\leq \sigma$, we deduce that $N/t_\cX(M)\in\cX$; and hence $N\in\cX$ which contradicts the choice of $t_\cX(M)$. Therefore, $M=t_\cX(M)\in\cX$. 
\end{proof}

\medskip
\begin{Proposition}
The category $\cA$ is semi-noetherian if and only if every non-zero object of $\cA$ contains a non-zero critical subobject.
\end{Proposition}
\begin{proof}
If $\cA$ is semi-noetherian, by \cref{excrit}, every non-zero module has a non-zero critical subobject. Conversely if $\cA_\tau$ is a proper localizing subcategory of $\cA$, then $\cA/\cA_\tau$ contains a non-zero object $Y$. Then  there exists a non-zero object $X$ in $\cA$ such that $t_\tau(X)=0$ and $F_\tau(X)=Y$ where $F_\tau:\cA\To\cA/\cA_\tau$ is the canonical functor. Hence, by the assumption, $X$ contains a non-zero $\lambda+1$-critical subobject $C$, where $\lambda\geq\tau$. Consequently, by \cref{critp}, we have $\cA_\tau\neq\cA_{\lambda+1}$ which is a contardiction.
\end{proof}

\section{The atom spectrum of a locally coherent  category}

 A finitely generated object $Y$ in $\cA$ is {\it finitely presented} if every epimorphism $f:X\To Y$ in $\cA$ with $X$ finitely generated has a finitely generated kernel $\Ker f$. A finitely presented object in $\cA$ is {\it coherent} if every its finitely generated subobject is finitely presented. We denote by fg-$\cA$, fp-$\cA$ and coh-$\cA$, the full subcategories of $\cA$ consisting of finitely generated, finitely presented and coherent objects, respectively. 

We recall that a Grothendieck category $\cA$ is {\it locally coherent} if every object in $\cA$ is a direct limit of coherent objects; or equivalently finitely generated subobjects of finitely presented objects are finitely presented. According to \cite[2]{Ro} and \cite{H} a Grothendieck category $\cA$ is locally coherent if and only if fp-$\cA=$coh-$\cA$ is an abelian category. A localizing subcategory $\cX$ of $\cA$ is said to be of {\it finite type} provide that the corresponding right adjoint  of the inclusion functor $\cX\To\cA$ commutes with direct limits.  
 If $\cA$ is a locally noetherian Grothendieck category, then fg-$\cA=$noeth-$\cA=$fp-$\cA=$coh-$\cA$ so that $\cA$ is locally coherent. In this case, any localizing subcategory $\cX$ of $\cA$ is of finite type. Throughout this section $\cA$ is a locally coherent category.
  
   In our investigation, we use  the following lemma due to Krause \cite[Lemma 1.1]{Kr}.
 
 \medskip

 \begin{Lemma}\label{lemkru}
 An object $X\in\cA$ is finitely generated if and only is for any epimorphism $\phi:Y\To X$ in $\cA$, there is a finitely generated subobject $U$ of $Y$ such that $\phi(U)=X$
 \end{Lemma}

Let $\alpha\in\ASpec\cA$. For any monoform objects $G,H$ in $\cA$ with $\alpha=\overline{H}=\overline{G}$, by the definition $H$ contains a non-zero subobject $K$ isomorphic to a non-zero subobject $L$ of $G$. As $G$ and $H$ are uniform objects,  $E(H)=E(K)=E(L)=E(G)$  is an indecomposable injective object of $\cA$. We denote $E(H)$ by $E(\alpha)$ as it is independent of the choice of a  monoform representative of $\alpha$. 

\medskip

\begin{Lemma}\label{llem}
Let $M$ be an object in $\cA$. Then  $\ASupp M=\{\alpha\in\ASpec\cA|\hspace{0.1cm} \Hom(M,E(\alpha))\neq 0\}$.
\end{Lemma}
\begin{proof}
 Given $\alpha\in\ASupp M$, there exist subobjects $K\subset L\subseteq M$ such that $H=L/K$ is a monoform object with $\alpha=\overline{H}$. Since $\Hom(H,E(\alpha))\neq 0$, we have $\Hom(L,E(\alpha))\neq 0$ and consequently $\Hom(M,E(\alpha))\neq 0$. For the other direction, let $\Hom(M,E(\alpha))$ be non-zero  with a non-zero morphism $f:M\To E(\alpha).$ Then $\alpha\in\ASupp \Im f\subseteq\ASupp M$ and so the inclusion holds.
\end{proof}

\medskip

For every subcategory $\cS$ of $\cA$, we denote by $\overset{\To}\cS$, the full subcategory of $\cA$ consisting of direct limits $\underset{\rightarrow}{\rm lim}X_i$ with $X_i\in\cS$ for each $i$. For $\cS\subset \fp\cA$, we denote by $\sqrt{\cS}$, the smallest Serre subcategory of $\fp\cA$ containing $\cS$. The following result due to Herzog [H] is crucial to our investigation.    
 
 \medskip
 \begin{Proposition}[\cite{H}, Proposition 2.3]\label{prh} Let $\cS$ be a Serre subcategory of $\fp\cA$ and $A$ is a finitely generated object of $\cA$. Then $A\in\overrightarrow{\cS}$ if and only if for every $B\in\fp\cA$, every epimorphism $\epsilon:B\To A$ factors through a quotient $S$ of $B$ which lies in $\cS$.
 \end{Proposition}
 
 The following lemma establishes an alternative topology on $\ASpec\cA$.
 \medskip
 \begin{Lemma}\label{bbas}
 The set $\{\ASupp M|\hspace{0.1cm} M\in\fp\cA\}$ forms a basis of open subsets for a topology on $\ASpec\cA$. Furthermore, if $M\in\fp\cA$, then $\ASupp M=\emptyset$ if and only $M=0$.
 \end{Lemma}
 \begin{proof}
 Since $\cA$ is locally coherent, it is clear that for every $\alpha\in\ASpec\cA$, there exists a finitely presented object $M$ in $\cA$ such that $\alpha\in\ASupp M$. If $M_1$ and $M_2$ are finitely presented objects in $\cA$ and $\alpha\in\ASupp M_1\cap\ASupp M_2$, then by \cref{llem}, there exist non-zero morphisms $f_i:M_i\To E(\alpha)$ for $i=1,2$. Since $E(\alpha)$ is uniform, $\Im f_1\cap\Im f_2$ is a non-zero subobject of $E(\alpha)$ and so it contains a non-zero finitely generated subobject $X$ as $\cA$ is locally coherent and so locally finitely generated. Therefore $X\in\overset{\longrightarrow}{\sqrt{M_1}}\cap \overset{\longrightarrow}{\sqrt{M_2}}$ and we have the following pull-back diagram
 $$\xymatrix{&&0\ar[d]&0\ar[d]\\
0\ar[r]&K\ar[r]\ar@{=}[d]& L\ar[d]\ar[r]^{\phi}& X\ar[r]\ar[d]\ar[r]& 0\\
0\ar[r]&K\ar[r]& M_1\ar[r]^{f_1}\ar[d]& \Im f_1\ar[r]\ar[d]& 0\\
 &&Z\ar@{=}[r]\ar[d]\ar@{=}& Z\ar[d]\\
&&0&0.}$$  
 
 By \cref{lemkru}, there exists a finitely generated subobject $L_1$ of $L$ such that $\phi(L_1)=X$. Since $M_1$ is finitely presented, $L_1$ is finitely presented as $\cA$ is locally coherent. By virtue of \cref{prh}, the morphism $\phi:L_1\To X$ factors  through a quotient $N$ of $L_1$ which lies in $\sqrt{M_2}$ as $X\in\overset{\longrightarrow}{\sqrt{M_2}}$. Therefore $N$ is finitely presented and $\alpha\in\ASupp N\subset\ASupp M_1\cap\ASupp M_2$. To prove the second claim, every non-zero finitely presented object $M$ in $\cA$ has a simple quotient object $S$ so that $\alpha=\overline{S}\in\ASupp  M$. 
 \end{proof}

 \medskip
 \begin{Remark}
 
 \begin{itemize}
\item To prevent any misunderstanding, we use the symbol $\ZASpec\cA$ instead of $\ASpec\cA$ with the new topology based on \cref{bbas}. To be more precise, if $\cU$ is an open subset 
of $\ZASpec\cA$, then there exists a family $(M_\alpha)$ in $\fp\cA$ such that $\cU=\bigcup_\alpha\ASupp M_\alpha$. For every subcategory $\cX$ of $\fp\cA$, we set $\ASupp\cX=\bigcup_{M\in\cX}\ASupp M$. 

\item We notice that m-$\ASpec\cA$ is a dense subset of $\ZASpec\cA$. Because if $M$ is a finitely presented object in $\cA$, it contains a maximal subobject $N$ so that the maximal atom $\overline{M/N}\in\ASupp M$.

\item We remark that if $M$ is a finitely presented object in $\cA$, then $\ASupp K$ is an open subset of $\ZASpec\cA$ for any subobject $K$ of $M$. Because $K=\bigcup_i K_i$ is the direct union of its finitely generated subobjects. Since $M$ is finitely presented, each $K_i$ is finitely presented so that $\ASupp K=\bigcup_i\ASupp K_i$ is an open subset
 of $\ZASpec\cA$.
  \end{itemize} 
\end{Remark}
  \medskip
  
\begin{Theorem}\label{corpropp}
The map $\cX\mapsto \ASupp\cX$  provides an inclusion-preserving bijective correspondence between Serre subcategories of $\fp\cA$ and  open subsets of $\ZASpec\cA$. The inverse map is given by $\cU\mapsto \ASupp^{-1}\cU\cap\fp\cA$.
\end{Theorem}
\begin{proof}
Assume that $\cX$ is a Serre subcategory of $\fp\cA$ and $\cU$ is an open subset of $\ZASpec\cA$. It is clear that $\ASupp\cX$ is an open subset of $\ZASpec\cA$ and $\ASupp^{-1}\cU\cap\fp\cA$ is a Serre subcategory of $\fp\cA$. Then the maps are well-defined. In order to prove $\ASupp^{-1}(\ASupp\cX)\cap\fp\cA=\cX$, it suffices to show that $\ASupp^{-1}(\ASupp\cX)\cap\fp\cA\subset\cX$ and the other direction is clear. Given $M\in \ASupp^{-1}(\ASupp\cX)\cap\fp\cA$, we have $\ASupp M\subset\ASupp\cX$ and so $\ASupp M\subset\ASupp\overset{\longrightarrow}{\cX}$.
 It follows from \cite[Proposition 2.15]{H} that $F(M)$ is a finitely presented object in $\cA/\overset{\longrightarrow}{\cX}$, where $F:\cA\To\cA/\overset{\longrightarrow}{\cX}$ is the canonical functor. Moreover, $\ASupp F(M)=\ASupp M\setminus \ASupp\cX=\emptyset$ by \cite[Proposition 5.6]{K}; and hence $F(M)=0$ by \cref{bbas}. Therefore, $M\in \overset{\longrightarrow}{\cX}$ and so $M=\underset{\rightarrow}{\rm lim}M_i$ is the direct limit of objects $M_i$ in $\cX$. Since $M$ is finitely presented, it is a direct summand of some $M_i$ so that $M\in\cX$. 
 
 Let $\alpha\in\ASupp(\ASupp^{-1}\cU\cap\fp\cA)$. Then there exists $M\in\ASupp^{-1}\cU\cap\fp\cA$  such that $\alpha\in\ASupp M\subset \cU$; consequently $\alpha\in\cU$. For the other direction, if $\alpha\in\cU$, then there exists $M\in\fp\cA$ such that $\alpha\in\ASupp M\subseteq\cU$. Hence $M\in\ASupp^{-1}\cU\cap\fp\cA$ so that $\alpha\in\ASupp(\ASupp^{-1}\cU\cap\fp\cA)$. Consequently $$\ASupp(\ASupp^{-1}\cU\cap\fp\cA)=\cU.$$
\end{proof}

\begin{Corollary}\label{ccooo}
The map $\cU\mapsto \overrightarrow{\ASupp^{-1}\cU\cap \fp\cA}$ provides an inclusion-preserving  bijective correspondence between open subsets of $\ZASpec\cA$ and  localizing subcategories of finite type of $\cA$. The inverse map is given by $\cX\mapsto \ASupp(\cX\cap\fp\cA)$.
\end{Corollary}
\begin{proof}
Given an open subset $\cU$ of $\ZASpec\cA$, by \cite[Theorem 2.8]{Kr}, $\overrightarrow{\ASupp^{-1}\cU\cap \fp\cA}$ is a localizing subcategory of finite type of $\cA$. It follows from \cref{corpropp} that $$\ASupp(\overrightarrow{\ASupp^{-1}\cU\cap \fp\cA}\cap\fp\cA)=\ASupp(\ASupp^{-1}\cU\cap \fp\cA)=\cU.$$ Also, for every localizing subcategory $\cX$ of finite type of $\cA$, it follows from \cref{corpropp} and \cite[Lemma 2.3]{Kr} that $$ \overrightarrow{\ASupp^{-1}(\ASupp(\cX\cap\fp\cA))\cap\fp\cA}=\overrightarrow{\cX\cap\fp\cA}=\cX.$$
\end{proof}

An open cover of a topological space $\cU$ is a family $\{\cU_\lambda|\hspace{0.1cm}\lambda\in\Lambda\}$ of open subsets of $\cU$ satisfying $\cU=\bigcup_{\lambda\in\Lambda}\cU_\lambda$. A topological space $\cU$ is {\it quasi-compact} if for any open cover $\{\cU_i\}_{\lambda\in\Lambda}$ of $\cU$, there exists a finite subset $\Gamma$ of $\Lambda$ such that $\cU=\bigcup_{\gamma\in\Gamma}\cU_\gamma$.

\medskip
\begin{Proposition}\label{qquasi}
 An open subset $\Phi$ of $\ZASpec\cA$ is a quasi-compact if and only if there exists an object $M$ in $\fp\cA$ such that $\Phi=\ASupp M$.
\end{Proposition}
\begin{proof}
Let $\Phi$ be quasi-compact. Since $\Phi$ is open, by definition, we have $$\Phi=\bigcup_{M\in\ASupp^{-1}\Phi\cap\fp\cA}\ASupp M.$$ Then there exists $M_1,\dots, M_n\in \ASupp^{-1}\Phi$ such that $\Phi=\bigcup_{i=1}^n\ASupp M_i=\ASupp(\bigoplus_{i=1}^n M_i)$. Hence the result follows as $\bigoplus_{i=1}^n M_i$ is finitely presented. Conversely, assume that $M$ is a finitely presented object in $\cA$ and $\{\Phi_\lambda|\hspace{0.1cm}\lambda\in\Lambda\}$ is an open cover of $\ASupp M$. It is clear that $$\bigcup_{\lambda\in\Lambda}\Phi_{\lambda}=\ASupp\sqrt{\langle\ASupp^{-1}\Phi_{\lambda}\cap\fp\cA|\hspace{0.1cm}\lambda\in\Lambda\rangle}$$ so that  $$\ASupp M=\bigcup_{\lambda\in\Lambda}\Phi_{\lambda}=\ASupp\sqrt{\langle\ASupp^{-1}\Phi_{\lambda}\cap\fp\cA|\hspace{0.1cm}\lambda\in\Lambda\rangle}.$$  Hence  \cref{corpropp} yields that $M\in\sqrt{\langle\ASupp^{-1}\Phi_{\lambda}\cap\fp\cA|\hspace{0.1cm}\lambda\in\Lambda\rangle}$. Then $M$ has a finite filtration of finitely presented subobjects $0=M_0\subset M_1\subset M_2\subset\dots\subset M_n=M$ such that $M_i/M_{i-1}$ is a subquotient of some finitely presented object $N_{\lambda_i}\in\ASupp^{-1}\Phi_{\lambda_i}$ for $1\leq i\leq n$. Therefore $\ASupp M=\bigcup_{i=1}^n\Phi_{\lambda_i}$   
\end{proof}

\medskip

\begin{Lemma}\label{mmon}
Any monoform object in $\fp\cA$ is a monoform object in $\cA$.
\end{Lemma}
\begin{proof}
 Let $H$ be a monoform object in $\fp\cA$. If $H$ is not monoform object in $\cA$, then there exists a subobject $N$ of $H$ and a  non-zero  subobject $Y$ of $H/N$ isomorphic to a subobject $X$ of $H$. Since $\cA$ is locally finitely generated, we may assume $X$ is finitely generated and since $H$ is finitely presented, $Y\cong X$ is finitely presented as $\cA$ is locally coherent. We now have the following pull-back diagram 
 
$$\xymatrix{&&0\ar[d]&0\ar[d]\\
0\ar[r]&N\ar[r]\ar@{=}[d]& K\ar[d]\ar[r]^{\phi}& Y\ar[r]\ar[d]\ar[r]& 0\\
0\ar[r]&N\ar[r]& H\ar[r]\ar[d]& H/N\ar[r]\ar[d]& 0\\
 &&Z\ar@{=}[r]\ar[d]\ar@{=}& Z\ar[d]\\
&&0&0.}$$  
 
Since  $Y$ is finitely generated, by \cref{lemkru},  there exists a finitely generated subobject $D$ of $K$ such that $\phi(D)=Y$. We observe that $D$ is a finitely generated subobject of $H$ and so is finitely presented. Set $C=\Ker(\phi\circ i)$ where $i:D\To K$ is the inclusion morphism. Consider a non-zero finitely generated subobject $N_1$ of $N=\Ker\phi$. Then $N_1$ and $D$ are non-zero finitely presented subobjects of $H$. Hence $N_1\cap D$ is a non-zero subobject of $C$ as $H$ is a monoform object in $\fp\cA$. Therefore $C$ is non-zero and since 
$\fp\cA$ is abelian, the exact sequence 
$0\To C\To D\To Y\to 0$ implies that $C$ is finitely presented; keeping in mind that $X\cong Y=D/C$. Now, this implies that $H$ and $H/C$ has a common non-zero finitely presented subobject $X$ which contradicts  the fact that $H$ is a monoform object of $\fp\cA$. 
\end{proof}

\medskip

\begin{Remark and Notation}
For a locally coherent  category $\cA$, the subcategory $\fp\cA$ is abelian. Thus, the atom spectrum of $\fp\cA$ can be investigated independently. If $\beta\in\ASpec\fp\cA$, then there exists a monoform object $H\in\fp\cA$ such that $\beta=\overline{H}$. It follows from \cref{mmon} that $H$ is a monoform object of $\cA$ so that $\alpha=\overline{H}\in\ASpec\cA$. To prevent any misunderstanding, we denote $\beta=\alpha_f$. By this notation, we will have
 \[\ASpec\fp\cA=\{\alpha_f=\overline{H}\mid\ \alpha=\overline{H}\in\ASpec\cA\text{ for some monoform object } H\in\fp\cA.\}\] 
For every object $M$ in $\fp\cA$, we use the symbol $\ASupp_f M$ for atom support of $M$ in $\ASpec\fp\cA$ instead of $\ASupp M$. Similarly we use the symbol $\AAss_f M$ instead of $\AAss M$. 

\end{Remark and Notation}

\medskip

\begin{Corollary}
 Let $M$ be a finitely presented object in $\cA$. Then \[\AAss_f M=(\AAss M)_f:=\{\alpha_f\in\ASpec\fp\cA\mid\  \alpha\in\AAss M\}.\]
 \end{Corollary}
 \begin{proof}
If $\alpha\in\AAss M$, then $M$ contains a monoform subobject $H$ in $\cA$ such that $\alpha=\overline{H}$. Since $\cA$ is locally coherent, we may assume that $H$ is finitely generated. Hence $H$ is finitely presented because $M$ is finitely presented. Moreover, it is clear that $H$ is a monoform object of $\fp\cA$ and consequently $\overline{H}=\alpha_f\in\AAss _fM$. To prove the other direction, if $\alpha_f\in\AAss_f M$, then $M$ contains a monoform subobject $H$ in $\fp\cA$ such that $\alpha_f=\overline{H}$. It follows from \cref{mmon} that $H$ is a monoform object in $\cA$ so that $\alpha\in\AAss M$. 
 \end{proof}

\medskip

\begin{Proposition}\label{homee}
 There is an injective continuous map $g:\ASpec\fp\cA\to\ZASpec\cA$ given by $\alpha_f\mapsto \alpha$.
  Moreover, if $\fp\cA$ is locally monoform, then  $\Im g$ is a dense subspace in $\ZASpec\cA$.
\end{Proposition}
\begin{proof}
Given $\alpha_f\in\ASpec\fp\cA$, there exists a monoform object $H\in\fp\cA$ such that $\alpha_f=\overline{H}$. Now, \cref{mmon} implies that $\alpha=\overline{H}\in\ZASpec\cA$ and so $g$ is well-defined. Moreover, it is clear that $g$ is injective and we show that $g$ is continuous. For any  finitely presented object $M$ in $\cA$, we prove that $g^{-1}(\ASupp M)=\ASupp_f M$.  The definition implies that $\ASupp_f M\subset g^{-1}(\ASupp M)$. To prove the converse, assume that $\alpha_f\in g^{-1}(\ASupp M)$. Then there exists a monoform object $H$ in $\fp\cA$ such that $\alpha=\overline{H}$ and $H$ is a subquotient of $M$. Using \cref{lemkru}, we can choose such $H$ such that $H=L/K$, where $L$ is a finitely presented subobject of $M$. Since $\fp\cA$ is abelian, we deduce that $K$ is finitely presented; and hence $\alpha_f\in\ASupp_f M$. Then $g^{-1}(\ASupp M)=\ASupp_f M$ and by \cref{remop}, $\{\ASupp_fM\mid M\in\fp\cA\}$ forms a basis for topology on $\ASpec\fp\cA$. Hence the induced topology by $g$ coincides with the topology on $\ASpec\fp\cA$. To prove the second claim, the assumption and \cref{mmon} imply that every finitely presented object $M$ contains a finitely presented monoform subobject $H$ in $\cA$. Setting $\overline{H}=\alpha_f$, we deduce that $\alpha=g(\alpha_f)\in\ASupp M\cap\Im\theta$.
\end{proof}

 It is known as Cohen's theorem that if prime ideals of a commutative ring $A$ are finitely generated, then $A$ is noetherian. In the following theorem we extend this theorem for locally coherent categories where we identify $\ASpec\fp\cA$ to the image of $g$ in \cref {homee}.

\begin{Theorem}\label{noeth}
If $\ASpec\fp\cA=\ASpec\cA$ (i.e. if every $\alpha\in\ASpec\cA$ has a monoform representative in $\fp\cA$), then $\cA$ is locally noetherian.
\end{Theorem}
\begin{proof}
Given a finitely presented object $M$ in $\cA$, it suffices to show that $M$ is noetherian. Otherwise, $M$ contains a subobject which is not finitely generated. It follows from \cite[Chap IV, Proposition 6.6]{St} that the class of subobjects of $M$ is a set. Then the set of subobjects of $M$ $$\{N\subset M|\hspace{0.1cm} N\hspace{0.1cm} is\hspace{0.1cm}{\rm not}\hspace{0.1cm} {\rm finitely}\hspace{0.1cm} {\rm generated} \}$$
 has a maximal element $K$ by Zorn's lemma. We claim that $M/K$ is a monoform object in $\cA$. Otherwise, there exists a subobject $L$ of $M$ containing $K$ such that $M/K$ and $M/L$ has a common non-zero subobject $X$. The maximality of $K$ implies that $L$ is finitely generated and so is finitely presented. Thus  $M/L$ is finitely presented so that we may assume that $X$ is finitely presented. Hence there exists a subobject $K_1$ of $M$ containig $K$ such that $X\cong K_1/K$ is finitely presented. Again the maximality implies that $K_1$ is finitely generated and so is finitely presented. Since $\fp\cA$ is abelian, we deduce that $K$ is finitely presented which is a contradiction. Therefore $M/K$ is monoform and so $\overline{M/K}\in\ASpec\cA$ and so the assumption implies that $M/K$ contains a finitely presented monoform subobject $Y/K$. Again by the maximality of $K$, the subobject $Y$ is finitely presented so that $K$ is finitely presented which is a contradiction.
\end{proof}
\medskip

 In the following theorem, we provide a sufficient condition under which a localizing subcategory of $\cA$ is of finite type.

\medskip
\begin{Theorem}\label{fgftg}
Let $\cX$ be a localizing subcategory of $\cA$ such that every atom in $\ASupp\cX$ has a  monoform representative in $\fp\cA$. Then $\cX$  is of finite type.
\end{Theorem}
\begin{proof}
By \cite[Theorem 2.8]{Kr}, it suffices to show that $\cX=\overrightarrow{\cX\cap\fp\cA}$. The inclusion $\overrightarrow{\cX\cap\fp\cA}\subset \cX$ is clear. For the other direction, let $M\in\cX$. Since $\cA$ is locally finitely generated, $M=\bigcup_iM_i$ is the direct union of its finitely generated subobjects. Then it suffices to show that each $M_i$ is finitely presented. Hence, we may assume that $M$ is a finitely generated object in $\cX$ and we have to prove that $M$ is finitely presented. Since $\cA$ is locally coherent, there exists a finitely presented object $N$ in $\cA$ and a subobject $L$ of $N$ such that $M=N/L$. Suppose that $L$ is not finitely generated. Then the set  
 $$\{X\subset N| \hspace{0.1cm}L\subset X \hspace{0.1cm}{\rm and}\hspace{0.1cm} X\hspace{0.1cm} is\hspace{0.1cm}{\rm not}\hspace{0.1cm} {\rm finitely}\hspace{0.1cm} {\rm generated} \}$$
  has a maximal element $K$ by Zorn's lemma. By a similar proof as given in \cref{noeth}, we can deduce that $N/K$ is a monoform object in $\cA$. Thus $\overline{N/K}\in\ASupp M$. Now the assumption and \cref{mmon} imply that $N/K$ contains a finitely presented subobject $Y$. Therefore by the maximality of $K$,  $N$ has a finitely presented subobject $Z$ containing $K$ such that $Y=Z/K$. Since $\fp\cA$ is abelian, we conclude that $K$ is finitely presented which is a contradiction. Hence $L$ is finitely generated so that $M$ is finitely presented.  
\end{proof}

\section{Ziegler spectrum of a locally coherent category}

Throughout this section $\cA$ is a locally coherent category. The collection of  the isomorphism classes of indecomposable injective objects in $\cA$, denoted by $\Zg\cA$, forms a set. This is because every indecomposable injective object is the injective envelope of some quotient of an element of the generating set of $\cA$.
  \begin{Definition}
 The {\it Ziegler support} $\ZSupp(M)$ of an object $M$ in $\cA$ is a subset of $\Zg\cA$, that is $$\ZSupp(M)=\{I\in\Zg\cA|\hspace{0.1cm} \Hom(M,I)\neq 0\}.$$  
\end{Definition}
For any exact sequence $0\To N\To M\To K\To 0$ of objects in $\cA$, we have $$\ZSupp(M)=\ZSupp (N)\cup\ZSupp(K).$$
For any subcategory $\cX$ of $\cA$, we define $\ZSupp(\cX)=\bigcup_{M\in\cX}\ZSupp(M)$ and for any subset $\cU$ of $\Zg\cA$, we define $\ZSupp^{-1}(\cU)=\{M\in\cA|\hspace{0.1cm} \ZSupp (M)\subset \cU\}$.

 Herzog \cite{H} proved that the collection $\{\cO(M)|\hspace{0.1cm} M\in\fp\cA\}$  satisfies the axioms for a basis of open subsets defining  a topology on $\Zg\cA$. The resulting  topological space is called the {\it Ziegler spectrum} of $\cA$. The assignment $\alpha\mapsto E(\alpha)$ makes  $\ZASpec\cA$ a subspace of $\Zg\cA$. More precisely we have the following lemma.

\medskip

\begin{Lemma}\label{ors}
There exists an injective continuous  map $\theta:\ZASpec\cA\To\Zg\cA$, given by $\alpha\mapsto E(\alpha)$. In particular, $\Im \theta$ is a dense subspace of $\Zg\cA$. 
\end{Lemma}
\begin{proof}
The injectivity of $\theta$ is clear. For every $M\in\fp\cA$, it follows from \cref{llem} that $\theta^{-1}(\cO(M))=\ASupp M$. Then $\theta^{-1}(\cO)$ is an open subset of $\ZASpec\cA$ for any open subset $\cO$ of $\Zg\cA$ and the induced topology by $\theta$ coincides with the topology on $\ZASpec\cA$. To prove the second claim, every finitely presented object $M$ has a maximal subobject $N$ so that $M/N$ is a monoform object in $\cA$. If we consider $\alpha=\overline{M/N}$, it is clear that $\Hom(M,E(\alpha))\neq 0$ and so $\theta(\alpha)\in\Im\theta\cap\cO(M)$. 
\end{proof}

A category $\cA$ is said to be {\it locally uniform} (or locally coirreducible by \cite[p. 330]{Po}) if every non-zero object in $\cA$ has a uniform subobject. The following proposition establishes a characterization of locally monoform objects.

\medskip

\begin{Proposition}\label{hemzig}
A category $\cA$ is locally monoform if and only if $\cA$ is locally uniform and the map $\theta:\ZASpec\cA\To\Zg\cA$, given by  $\alpha\mapsto E(\alpha)$ is a homeomorphism. 
\end{Proposition}
\begin{proof}
Let $\cA$ be locally monoform. Then it is clear that $\cA$ is locally uniform. Let  $E$ be an  indecomposable injective object in $\cA$. By the assumption, the object $E$ contains a monoform subobject $H$ and so $E=E(\alpha)$, where $\alpha=\overline{H}$. This implies that $\theta$ is surjective. Therefore, it follows from \cref{ors} that $\theta(\ASupp M)=\cO(M)$ for every $M\in\fp\cA$ so that $\theta$ is a homeomorphism.  Conversely, let $M$ be an object in $\cA$. Since $\cA$ is locally uniform, $M$ contains a uniform subobject $U$ so that $E(U)$ is an injective indecomposable. Since $\theta$ is surjective, there exists a monoform object $\alpha\in\ASpec\cA$ such that $E(\alpha)=E(U)$. Then there exists a  monoform object $H$ in $\cA$ such that $\alpha=\overline{H}$ and $E(H)=E(U)$ so that $H\cap U$ is a monoform subobject of $M$.  
\end{proof}

For a commutative ring $A$,  We have the following easy lemma.

\begin{Lemma}\label{llkk}
 Let $\frak p$ be an ideal of $A$. Then $A/\frak p$ is a monoform $A$-module if and only if $\frak p$ is a prime ideal of $A$. 
\end{Lemma}
\begin{proof}
Clear.
\end{proof}

We recall that a topological space $X$ is said to be {\it $T_0$-space} (or {\it Kolmogorov}) if for any distinct points $x,y$ of $X$, there exists an open subset of $X$ containing exactly one of them. The following example \cite{GP}, pointed out by T. Kucera, shows that a locally coherent category $\cA$ need not be locally monoform and also $\Zg\cA$ need not be $T_0$-space.

\begin{Example}\label{eex}
Let $A=k[X_n(n\in\omega)]$ be a polynomial ring over a field $k$ in infinitely many commuting indeterminates. Then $A$ is clearly a commutative coherent ring. Let $I=(X_n^{n+1}: n\in\omega)$. Then $I$ is not prime as $X_1\notin I$ but $X_1^2\in I$. According to \cite{P}, the injective $A$-module $E(A/I)$ is indecomposable and it does not have the form $E(A/\frak p)$ for any prime ideal $\frak p$ of $A$ by \cite[9.1]{P}. In \cref{hemzig}, if $f$ is surjective, then there exists $\alpha\in\ASpec\Mod A$ such that $E(A/I)=E(\alpha)$. Hence there exists a cyclic monoform $A$-module $H$ such that $E(H)=E(A/I)$. But \cref{llkk} implies that $H=A/\frak p$ for some prime ideal $\frak p$ of $A$ which is a contradiction.  Therefore $\Mod A$ is not locally monoform by \cref{hemzig}. Furthermore, since $E(A/I)$ does not have the form $E(A/\frak p)$ for any prime ideal, by \cite[Theorem 1.4]{GP}, there exists a prime ideal $\frak p$ of $A$ such that $\Lambda(E(A/I))=\Lambda(E(A/\frak p))$, where for any injective indecomposable $E$, the set $\Lambda(E)$ denotes the intersection of all open subsets of $\Zg\Mod A$ containing $E$. But this implies that $\Zg \Mod A$ is not $T_0$-space. 
\end{Example} 

 The following proposition provides a sufficient condition under which a locally coherent category is semi-noetherian (and so locally monoform).
\begin{Proposition}\label{semii}
If $\cA_\tau$ is a localizing subcategory of finite type of $\cA$, then $\cA$ is semi-noetherian. 
\end{Proposition}
\begin{proof}
If $\cA$ is not semi noetherian, then $\cA_\tau$ is a proper subcategory of $\cA$. It follows from \cite[Theorem 2.16]{H} that $\cA/\cA_\tau$ is locally coherent so that it contains a simple object $S$. Hence there exists an object $X\in\cA_{\tau+1}\setminus\cA_\tau $ such that $F_\tau(X)=S$ which is a contradiction. 
\end{proof}

\medskip
 Let $\Sp\cA$ be the set of  isomorphism classes of indecomposable injective objects in $\cA$. Krause \cite{Kr} has constructed  a topology on $\Sp\cA$ in which  for a subset $\cU$ of $\Sp\cA$, the closure of $\cU$ is defined as $\overline{\cU}=\left\langle^\bot \cU\cap \fp\cA\right\rangle^\bot$. The subsets $\cU$ of $\Sp\cA$ satisfying $\cU=\overline{\cU}$ form the closed subsets of $\Sp\cA.$ The following proposition shows that the Ziegler topology and the topology defined by Krause are identical.

\medskip

\begin{Lemma}\label{alft}
The Ziegler topology and Krause topology on $\Zg(\cA)$ are the same. 
\end{Lemma}
\begin{proof}
We show that $\Zg\cA$ and $\Sp\cA$ have the same open subsets. Given an open subset $\cO$ of $\Zg(\cA)$, it suffices to show that  $\left\langle^\bot\cO^c\cap\fp\cA\right\rangle^\bot=\cO^c$ and so $\cO^c$ will be a closed subset of $\Sp\cA$, where $\cO^c=\Sp\cA\setminus\cO$. If $I\in\cO^c$, then it is clear that $\Hom(^\bot\cO^c\cap\fp\cA,I)=0$ and so
 $I\in\left\langle^\bot\cO^c\cap\fp\cA
\right\rangle^\bot$. Conversely, if $I\in \left\langle^\bot\cO^c\cap\fp\cA\right\rangle^\bot\setminus\cO^c$, there exists $M\in\fp\cA$ such that $I\in\cO(M)\subseteq\cO$; and hence $M\notin^\bot\cO^c\cap\fp\cA$. Then $\Hom(M,\cO^c)\neq 0$ so that there exists $J\in\cO^c$ such that $\Hom(M,J)\neq 0$. But this implies that $J\in\cO(M)\subseteq \cO$ which is a contradiction. Now suppose that $\cO$ is an open subset of $\Sp\cA$ and so $\cO^c=\left\langle^\bot\cO^c \cap\fp\cA\right\rangle^\bot$. We now show that $\cO$ is an open subset of $\Zg\cA$. Given $I\in\cO$, we have $\Hom(^\bot\cO^c \cap\fp\cA,I)\neq 0$ and so  there exists $M\in^\bot\cO^c \cap\fp\cA$  such that $\Hom(M,I)\neq 0$. Thus $I\in\cO(M)$ and $\Hom(M,\cO^c)=0$. For every $J\in\cO(M)$, we have $\Hom(M,J)\neq 0$ which implies that $J\in\cO$. Therefore, $\cO(M)\subset \cO$ ; and consequently $\cO$ is an open subset of $\Zg\cA$.     
\end{proof}

Although the following proposition has been  investigated indirectly by Herzog [H] and Karuse [Kr], we give a proof because it will be used in the next section.

\begin{Proposition}\label{cortype}
The map $\cX\mapsto\ZSupp(\cX)$ provides an inclusion-preserving bijective correspondence between localizing subcategories of finite type of $\cA$ and  open subsets of $\Zg\cA$. The inverse map is $\cU\mapsto \ZSupp^{-1}(\cU)$.   
\end{Proposition}
\begin{proof}
Given a localizing subcategory $\cX$ of finite type of $\cA$, by \cite[Lemma 2.3]{Kr}, we have $\cX=\overset{\To}\cS$, where $\cS=\cX\cap\fp\cA$. Then $\ZSupp(\cX)=\bigcup_{M\in\cS}\ZSupp(M)$ that is an open subset of $\Zg\cA$. Given an open subset $\cU$ of $\Zg\cA$, it is clear by the definition that $\ZSupp^{-1}(\cU)=^\bot\cU^c$. Then \cite[Corollary 4.3]{Kr} and \cref{alft} imply that  $\ZSupp^{-1}(\cU)$ is a localizing subcategory of finite type of $\cA$.  To prove $\cX=\ZSupp^{-1}(\ZSupp(\cX))$, clearly $\cX\subset \ZSupp^{-1}(\ZSupp(\cX))$. For the other side, by the previous argument, $\ZSupp^{-1}(\ZSupp(\cX))$ is a localizing subcategory of finite type of $\cA$. Thus, for every $M\in \ZSupp^{-1}(\ZSupp(\cX))$, we have $M=\underset{\To}{\rm lim}M_i$ where each $M_i$ belongs to $\ZSupp^{-1}(\ZSupp(\cX))\cap \fp\cA$ by \cite[Lemma 2.3]{Kr}. Then $\ZSupp(M_i)\subset\ZSupp(\cX)$ for each $i$. Fixing $i$, since $\ZSupp (M_i)$ is a quasi-compact open subset of $\Zg\cA$ by \cite[Corollary 4.6]{Kr}, there exists $N\in\cX\cap \fp\cA$ such that $\ZSupp(M_i)\subset \ZSupp(N)$ and hence it follows from [H, Corollary 3.12] that $M_i\in\sqrt{N}$, where $\sqrt{N}$ is the  smallest Serre subcategory of $\fp\cA$ containing $N$. Clearly $\sqrt{N}\subset\cX$ and hence $M_i\in\cX$. Finally, this forces that $M\in\cX$ as $M=\underset{\To}{\rm lim}M_i$. 

It is clear that $\ZSupp(\ZSupp^{-1}(\cU))\subset\cU$. On the other hand, for every $I\in\cU$, there exists a finitely presented object $M$ such that $I\in\cO(M)\subset \cU$. This implies that $M\in\ZSupp^{-1}(\cU)$ and so the previous argument forces that $\ZSupp (M)\subset \ZSupp(\ZSupp^{-1}(\cU))$. Therefore $I\in \ZSupp(\ZSupp^{-1}(\cU))$; and hence $\cU=\ZSupp(\ZSupp^{-1}(\cU))$.    
\end{proof}

\section{In the case of commutative coherent rings}

In this section, we assume that $A$ is a commutative  coherent ring and we denote the category of $A$-modules by $\Mod A$. Also, we denote $\ZASpec \Mod A$, $\Zg\Mod A$ and $\fp\Mod A$ by $\ZASpec A$, $\Zg A$ and $\fp A$, respectively. 

\medskip

We recall that a subspace $Y$ of a topological space $X$ is {\it retract} provided that there exists a continuous function $r:X\to Y$ such that $r(y)=y$ for all $y$ in $Y$. Let $E$ be any indecomposable injective $A$-module. Set $\frak p_E$ to be the sum of annihilator ideals of non-zero elements of $E$. As $E$ is uniform, in view of \cite[9.2]{P}, it is easy to check that $\frak p_E$ is a prime ideal of $A$ (also one can see the paragraph before \cite[Theorem 1.4]{GP}). If we identify every $\alpha\in\ZASpec A$ with $E(\alpha)$ in $\Zg A$, then we have the following result.

\medskip
\begin{Proposition}
$\ZASpec A$ is a retract subspace of $\Zg A$.  
\end{Proposition}
\begin{proof}
We define $\theta:\Zg A\To \ZASpec A$ given by $E\To E(A/{\frak p_E})$. If $E(\alpha)\in\ZASpec\cA$, then by the above argument, there exists a prime ideal $\frak p=\frak p_{E(\alpha)}$ of $A$ such that $\alpha=\overline{A/\frak p}$ by \cref{llkk}. Hence $E(\alpha)=E(A/\frak p)$ so that $\theta(E(\alpha))=E(\alpha)$. For every finitely presented $A$-module $M$, we show that $\theta^{-1}(\ASupp M)=\cO(M)$. If $E\in\theta^{-1}(\ASupp M)$, we have $\theta(E)=E(A/\frak p_E)\in\ASupp M$ so that $E(A/\frak p_E)\in\cO(M)$ by \cref{llem}. It now follows from \cite[Theorem 1.4]{GP} that $E\in\cO(M)$. Conversely if $E\in\cO(M)$, then  $\theta(E)=E(A/\frak p_E
)\in\ASupp M$ by [GP, Theorem 1.4]. Therefore $E\in\theta^{-1}(\ASupp M)$.  
\end{proof}

\medskip

\begin{Lemma}\label{prb}
The set $\cF\cI=\{\ASupp A/\frak a|\hspace{0.1cm} \frak a$ is a finitely generated ideal of $A\}$ forms a basis for the topology 
 $\ZASpec A$.
\end{Lemma}
\begin{proof}
We observe that $\ASpec A=\ASupp A$. For any  finitely generated ideals $\frak a$ and $\frak b$ of $A$ and $\alpha\in\ASupp A/\frak a\cap\ASupp A/\frak b$. By \cref{llkk}, there exists a prime ideal $\frak p$ such that $\overline{A/\frak p}=\alpha$. It follows from \cref{llem} that $\Hom_A(A/\frak a,E(A/\frak p))\neq 0$ and $\Hom_A(A/\frak b,E(A/\frak p))\neq 0$. This implies that $\frak p\in V(\frak a)\cap V(\frak b)=V(\frak a+\frak b)$; and hence we have  $\alpha\in\ASupp A/(\frak a+\frak b)$. Therefore   
$\ASupp A/\frak a\cap\ASupp A/\frak b\subseteq\ASupp A/(\frak a+\frak b)$ and an easy argument implies that $\ASupp A/\frak a\cap\ASupp A/\frak b=\ASupp A/(\frak a+\frak b)$. On the other hand, assume that $M$ is a finitely presented $A$-module. Then $M$ is finitely generated, say $M=<x_1,\dots,x_n>$.  Set $M_0=0$ and $M_k=<x_1,\dots x_k>$ for $1\leq k\leq n$. Then $M_i/M_{i-1}$ is cyclic and finitely presented for each $i$. Thus there exist finitely generated ideals $\frak a_i$ such that $M_i/M_{i-1}\cong A/\frak a_i$ for all $1\leq i\leq n$. It is clear that $\ASupp M=\bigcup_{i=1}^n\ASupp A/\frak a_i$. 
\end{proof}

\medskip
For a spectral topological space $X$, Hochster \cite{Ho} endowed the underlying set with a new, dual topological by defining its open subsets as those of the form $Y=\bigcup_{i\in\Omega} Y_i$, where $X\setminus Y_i$ is a quasi-compact open subset of $X$ for each $i\in\Omega$. The symbol $X^\star$ denotes $X$  with the new topology. We write $\Spec^\star A$ for $(\Spec A)^\star$.

\medskip
\begin{Proposition}\label{pla}
The maps  $$\Spec^\star A\supseteq\cV\stackrel{\cO_{(-)}}\mapsto \cO_{\cV}=\{\alpha=\overline{A/\frak p}\in\ZASpec A|\hspace{0.1cm} \frak p\in\cV\}\subseteq \ZASpec A\hspace{0.15cm} and$$$$\ZASpec A\supseteq\hspace{0.1cm}\cO\stackrel{\cV_{(-)}}\mapsto \cV_{\cO}=\{\frak p\in\Spec A|\hspace{0.1cm}\overline{A/\frak p}\in\cO\}\subseteq \Spec^\star A$$
establish an inclusion-preserving  bijective correspondence between open subsets of $\Spec^\star A$ and those of $\ZASpec A$. 
\end{Proposition}
\begin{proof}
For every finitely generated ideal $\frak a$ of $A$, it is known that $\Spec A\setminus V(\frak a)$ is a qausi-compact open subset of $\Spec A$. Therefore for every open subset $\cV$ of $\Spec^\star A$, we have
$\cV=\bigcup_{\gamma\in\Gamma} V(\frak a_\gamma)$ where $\frak a_\gamma$ ranges over finitely generated ideals of $A$. By the same argument used in the proof of \cref{prb}, we have $\cO_{V(\frak a_\gamma)}=\ASupp A/\frak a_\gamma$ for all $\gamma\in\Gamma$; and hence \[\cO_\cV=\bigcup\cO_{V(\frak a_\lambda)}=\{\overline{A/\frak p}\in\ZASpec A|\hspace{0.1cm} \frak p\in\cV\}=\bigcup_{\gamma\in\Gamma}\ASupp A/\frak a_\gamma\] is an open subset of $\ZASpec A$. Therefore, the map $\cO_{(-)}$ is well-defined. Now, assume that  $\cO$ is an open subset of $\ZASpec\cA$. By \cref{prb}, we have  $\cO=\bigcup_{\Lambda}\ASupp A/\frak a_{\lambda}$, where $\frak a_\lambda$ ranges over finitely generated ideals.  Hence $\{\frak p\in\Spec A|\hspace{0.1cm}\overline{A/\frak p}\in\cO\}=\bigcup_{\lambda\in\Lambda} V(\frak a_\lambda)$ is a open subset of $\Spec^\star A$; and therefore the map $\cV_{(-)}$ is well-defined. Now, it is easy to see that $\cO_{\cV_{\cO}}=\cO$ and $\cV_{\cO_{\cV}}=\cV$.
\end{proof}

\medskip

\begin{Lemma}\label{prbr}
 The set $\{\ZSupp(A/\frak a)|\hspace{0.1cm} \frak a$ is a finitely generated ideal of $\cA\}$ forms a basis for the topology 
 $\Zg A$.
\end{Lemma}
\begin{proof}
We show that $\Zg A=\bigcup_{\frak a\in\text{fg-}A} \cO(A/\frak a)$. Any indecomposable injective module $E$ contains a cyclic submodule $A/\frak c$. This implies that  $E\in\ZSupp(A/\frak c)$. If $\frak c=0$, there in nothing to prove. Otherwise, $\frak c$ contains a finitely generated ideal $\frak d$ so that $E\in\cO(A/\frak d)$. The rest of the proof is similar to that of \cref{prb}.
\end{proof}

\medskip
A set $\kF$ of ideals of $A$ is said to be a {\it Gabriel topology} on $A$ if it satisfies:

(1). If $\frak a\in\kF$ and $\frak a\subset \frak b$, then $\frak b\in\kF$.

(2). If $\frak a$, $\frak b$ belong to $\kF$, then $\frak a\cap \frak b\in\kF$.

(3). If $\frak a\in\kF$ and $a\in A$, then $(\frak a:a)\in\kF$.

(4). If $\frak a$ is an ideal of $A$ and there exists $\frak b\in\kF$ such that $(\frak a:b)\in\kF$ for every $b\in\frak b$, then $\frak a\in \kF$. 

\medskip
A Gabriel topology $\kF$ is said to have a basis of finitely generated ideals, if  any ideal $\frak a\in\kF$ contains a finitly generated ideal $\frak b\in\kF$. 
Let $E$ be an indecomposable injective module in $\Zg A$. Clearly, the Gabriel topology $$\kF_E=\{\frak a|\hspace{0.1cm} E\notin\ZSupp(A/\frak a)\}$$  is associated to $\{M\in\Mod A|\hspace{0.1cm}\Hom_A(M,E)=0\}$, the localizing subcategory cogenerated by $E$.
More generally, let $\cE$ be a subset of $\Zg A$. To $\cE$, we associated a Gabriel topology $$\kF_\cE=\bigcap_{E\in\cE}\kF_E=\{\frak a|\hspace{0.1cm} \ZSupp(A/\frak a)\cap\cE=\emptyset\}.$$ On the other hand, for every Gabriel topology $\kF$ on $A$, we associate
$$D(\kF)=\{E\in\Zg A|\hspace{0.1cm} E\notin\ZSupp (A/\frak a)\hspace{0.1cm} {\rm for \hspace{0.1cm} every} \hspace{0.1cm}\frak a\in\kF\}=\Zg A\setminus\bigcup_{\frak a\in\kF}\ZSupp(A/\frak a).$$
Given a Gabriel topology $\kF$, an $A$-module $M$ is said to be $\kF$-{\it discrete} if $\Ann(x)\in\kF$ for all $x\in M$. The following theorem provides a new classification of Gabriel topologies on $A$.  

\begin{Lemma}\label{fgg}
Let $D$ be a closed subset of $\Zg A$. Then $\kF_D$ has  a basis of finitely generated ideals of $A$.
\end{Lemma}
\begin{proof}
We observe that $\Zg A\setminus D$ is an open subset of $\Zg A$; and hence by \cref{prbr}, there exists a family of finitely generated ideals $\{\frak b_{\lambda}|\hspace{0.1cm}\lambda\in \Lambda\}\subseteq \kF_D$ such that $\Zg A\setminus D=\bigcup_{\lambda\in\Lambda}\cO (A/\frak b_\lambda)$. Now, given an ideal $\frak a \in\kF_D$, by the definition, $\cO (A/\frak a)\subseteq \bigcup_{\lambda\in \Lambda}\cO (A/\frak b_\lambda)$ so that $A/\frak a$ is a quotient of $A/\frak c\in\sqrt{\langle A/\frak b_\lambda|\hspace{0.1cm}\lambda\in \Lambda\rangle}$ by \cref{cortype}. This implies that $\frak c\subseteq\frak a$ and  there exists a finite filtration of submodules of $A/\frak c$ 
$$0=X_0\subset X_1\subset\dots\subset X_{n-1}\subset X_n=A/\frak c$$ such that  each $X_i/X_{i-1}$ is a subquotient of $A/\frak b_{\lambda_i}$ for some $\lambda_i\in \Lambda$. Setting $\frak b=\frak b_{\lambda_1}\frak b_{\lambda_2}\dots \frak b_{\lambda_n}$, it follows that $\frak b\subset\frak c\subset\frak a$ and $\frak b\in\kF_D$ by \cite[Chap VI, p. 147, Lemma 5.3]{St}. 
\end{proof}

There are many classifications of Gabriel topologies on $A$ (e.g., see \cite{G, St}). For commutative coherent rings, the following theorem provides a new classification of Gabriel topologies on $A$. The subcategory of finitely generated modules is denoted by fg-$A$. 
 \medskip

\begin{Theorem}\label{flfg}
\begin{enumerate}
\item For any finitely generated ideal $\frak a$, the assignments $$\ASupp A/\frak a\longleftrightarrow\Zg A\setminus \cO(A/\frak a)$$ induce an inclusion-reversing bijection between open subsets of $\ZASpec A$ and closed subsets of $\Zg A$.

\item  The assignment $D\stackrel{\kF}\mapsto \kF_D$ induces  an inclusion-reversing bijective  between closed subsets of $\Zg A$ and Gabriel topologies having bases of finitely generated ideals of $A$. The inverse map is given by $\kF\stackrel{D}\mapsto D(\kF)$.

\item  The assignment   $\kF\mapsto \cX_\kF=\overrightarrow{\sqrt{\langle A/\frak a|\hspace{0.1cm} \frak a\in \kF\cap\fg A \rangle}}$ induces  an inclusion-preserving bijective  between  Gabriel topologies having bases of finitely generated ideals of $A$ and localizing subcategories of finite type of $\Mod A$. The inverse map is given by $\cX\mapsto\kF_{\cX}=\{\frak a|\hspace{0.1cm}A/\frak a\in\cX\}$.
\end{enumerate}
 
 \end{Theorem}
 \begin{proof}
 (1): For open subsets $\cU$ and $\cV$ of $\ZASpec\cA$ with $\cU\subseteq \cV$ there exist the families of finitely generated ideals $\{\frak b_{\lambda}|\hspace{0.1cm}\lambda\in \Lambda\}$ and $\{\frak a_{\gamma}|\hspace{0.1cm}\gamma\in \Gamma\}$ such that $\cU=\bigcup_{\lambda\in\Lambda}\ASupp(A/\frak b_\lambda)$ and $\cV=\bigcup_{\gamma\in\Gamma}\ASupp(A/\frak a_\gamma)$. For each $\lambda\in\Lambda$, the inclusion $\cU\subseteq\cV$ and  \cref{corpropp} imply that $A/\frak b_\lambda\in\sqrt{\langle\frak a_\gamma|\gamma\in\Gamma\rangle}$ so that $\cO(A/\frak b_\lambda)\subseteq\bigcup_{\gamma\in\Gamma}\cO(A/\frak a_\gamma)$. Therefore $\bigcup_{\lambda\in\Lambda}\cO(A/\frak b_\lambda)\subseteq\bigcup_{\gamma\in\Gamma}\cO(A/\frak a_\gamma)$ so that the first map reverses  inclusion. By a similar argument and using \cref{cortype}, the inverse map reverses the inclusion-reversion. Now the rest of claim is clear  \cref{ccooo}, \cref{cortype},
  \cref{prb} and \cref{prbr}.
  
 (2): It is clear that the maps are inclusion-reversing. Let $\kF$ be a Gabriel topology on $A$ having a basis of finitely generated ideal of $A$ and  $D$ be a closed subset of $\Zg A$. Then \cref{fgg} implies that $\kF_D$ has a basis of finitely generated ideals and $D(\kF)=\Zg A\setminus\bigcup_{\frak a\in\kF}\ZSupp(A/\frak a)=\Zg A\setminus\bigcup_{\frak a\in\kF\cap\fg A}\ZSupp(A/\frak a)$ is a closed subset of $\Zg A$ by \cref{prbr}. 
 
  We show that $\kF_{D(\kF)}=\kF$. Given $\frak a\in\kF$, the definition implies that $\Hom_A(A/\frak a,D(\kF))=0$ so that $\frak a\in\kF_{D(\kF)}$. Conversely, if $\frak b\in\kF_{D(\kF)}$, then $\Hom_A(A/\frak b,D(\kF))=0$ so that $\ZSupp(A/\frak b)\subseteq \bigcup_{\frak a\in\kF}\ZSupp (A/\frak a)$. It then follows from \cref{fgg} that  $\ZSupp(A/\frak b)\subseteq \bigcup_{\frak a\in\kF\cap{\text{fg-}}A}\ZSupp (A/\frak a)=\cO(\overrightarrow{\sqrt{\langle A/\frak a|\hspace{0.1cm} \frak a\in \kF\cap{\text{fg-}}A \rangle}})$. Now, \cref{cortype} implies that $A/\frak b\in\overrightarrow{\sqrt{\langle A/\frak a|\hspace{0.1cm} \frak a\in \kF\cap{\text{fg-}}A \rangle}}$. Therefore,  $A/\frak b$ is a quotient of some $A/\frak a$ for $\frak a\in\kF\cap{\text{fg-}}A$.  This implies that $\frak a\subseteq \frak b$ so that $\frak b\in\kF$.
 
  To prove $D(\kF_D)=D$, since $\Zg A\setminus D$ is an open subset of $\Zg\cA$, by \cref{prbr}, there exists a set of finitely generated ideals $\{\frak a_\lambda|\hspace{0.1cm}\lambda\in\Lambda\}$ of $A$ such that $\Zg A\setminus D=\bigcup_{\lambda\in\Lambda}\ZSupp(A/\frak a_\lambda)$. The inclusion     $D\subseteq D(\kF_D)$ is clear. Conversely, for any $E\in D(\kF_D)$, since $\frak a_\lambda\in\kF_D$ for any $\lambda\in\Lambda$,  we have $E\notin\ZSupp(A/\frak a_\lambda)$ for any $\lambda\in\Lambda$. Consequently $E\in D$.
 
 (3): It is clear that the maps are inclusion-preserving. Given a Gabriel topology $\kF$ with a basis of finitely generated ideals of $A$, it follows from \cite[Theorem 2.8]{Kr} that $\cX_{\kF}$ is a localizing subcategory of finite type of $\cA$. Moreover, for a localizing subcategory $\cX$ of $\Mod A$ of finite type, it follows from \cite[Chap VI, p. 146,  Theorem 5.1]{St} that $\kF_{\cX}$ is a Gabriel topology on $A$. If $\frak a\in\kF_\cX$, then we have $A/\frak a=\underset{\rightarrow}{\rm lim}M_i$ is the direct limit of objects $M_i$ in $\cX\cap\fp A$. Then $\frak a$ contains a finitely generated ideal $\frak b$  such that $A/\frak b\in\cX$. Hence $\kF_{\cX}$ has a basis of finitely generated ideals of $A$.
 
  We first show that $\kF_{\cX_\kF}=\kF$. Given $\frak a\in\kF$, it contains a finitely generated ideal $\frak b\in\kF$ so that $A/\frak b\in\cX$. Then $A/\frak a\in\cX_\kF$; and hence $\frak a\in\kF_{\cX_\kF}$. Conversely, if $\frak a  \in\kF_{\cX_\kF}$, then $A/\frak a\in\cX_\kF$. By a similar argument as mentioned above, $\frak a$ contains a finitely generated ideal $\frak b\in\kF$ so that $\frak a\in\kF$.
  
   Now, we show that $\cX_{\kF_\cX}=\cX$. Given $M\in\cX$, since $\cX$ and $\cX_{\kF_\cX}$ are localizing subcategory of finite type of $\cA$, we may assume that $M$ is finitely presented. Assume that  $M=<x_1,\dots,x_n>$, $M_0=0$ and $M_k=<x_1,\dots x_k>$ for $1\leq k\leq n$. Then $M_i/M_{i-1}=A/\frak a_i\in\cX$ for $1\leq i\leq n$. Hence $M_i/M_{i-1}\in\cX_{\kF_\cX}$ for $1\leq i\leq n$; and consequently $M\in\cX_{\kF_\cX}$. 
 \end{proof}







\end{document}